\documentclass[11 pt]{article}
\usepackage{geometry}
\usepackage{fullpage}
\usepackage{graphicx}
\usepackage{amssymb}
\usepackage{epstopdf}
\usepackage{amsmath,amsthm,amscd,amssymb}
\usepackage{latexsym}
\usepackage[colorlinks,citecolor=red,pagebackref,hypertexnames=false]{hyperref}

\numberwithin{equation}{section}

\theoremstyle{plain}
\newtheorem{theorem}{Theorem}[section]
\newtheorem{lemma}[theorem]{Lemma}
\newtheorem{corollary}[theorem]{Corollary}

\theoremstyle{definition}

\theoremstyle{remark}
\newtheorem{remark}[theorem]{Remark}

\newtheorem{case[theorem]}{Case}

\title{Maximal operators: scales, curvature and the fractal dimension }
\author{A. Iosevich, B. Krause, E. Sawyer, K. Taylor and I. Uriarte-Tuero}

\begin{document}
\maketitle

\begin{abstract} We establish $L^p$ bounds for the Bourgain-Stein spherical maximal operator in the setting of compactly supported Borel measures $\mu, \nu$ satisfying natural local size assumptions $\mu(B(x,r)) \leq Cr^{s_{\mu}}, \nu(B(x,r)) \leq Cr^{s_{\nu}}$. Taking the supremum over all $t>0$ is not in general possible for reasons that are fundamental to the fractal setting, but we are able to obtain single scale ($t \in [1,2]$) results. The range of possible $L^p$ exponents is, in general, a bounded open interval where the upper endpoint is closely tied with the local smoothing estimates for Fourier Integral Operators. 

In the process, we establish $L^2(\mu) \to L^2(\nu)$ bounds for the convolution operator $T_{\lambda}f(x)=\lambda*(f\mu)$, where $\lambda$ is a tempered distribution satisfying a suitable Fourier decay condition. More generally we establish a transference mechanism which yields $L^p(\mu) \to L^p(\nu)$ bounds for a large class of operators satisfying suitable $L^p$-Sobolev bounds. This allows us to effectively study the dimension of a blowup set ($\{x: Tf(x)=\infty \}$) for a wide class of operators, including the solution operator for the classical wave equation. Some of the results established in this paper have already been used to study a variety of Falconer type problems in geometric measure theory.  \end{abstract} 

\section{Introduction}

\vskip.125in 

The spherical maximal operator, 
$$ {\mathcal A}f(x)=\sup_{t>0} |A_tf(x)|=\sup_{t>0} \left| \int f(x-ty) d\sigma(y) \right|,$$ where $\sigma$ is the surface measure on the unit sphere, is a classical object that appears in variety of contexts in harmonic analysis, geometric measure theory, partial differential equation and geometric combinatorics. See, for example, \cite{St93}, \cite{So93}, \cite{Falc86} and \cite{IJL09}. It is known that 
\begin{equation} \label{bourgstein} {\mathcal A}: L^p({\Bbb R}^d) \to L^p({\Bbb R}^d) \ \text{for} \ p>\frac{d}{d-1}. \end{equation} 

This fact was established by Stein in dimension three and higher (\cite{St76}), and by Bourgain in dimension two (\cite{B86}). 

One of the motivations behind the Bourgain-Stein spherical maximal theorem is the classical initial value problem for the wave equation. Consider the equation 
\begin{equation} \label{wave} \Delta u=\frac{\partial^2 u}{\partial t^2}; \ u(x,0)=0; \ \frac{\partial u}{\partial t}(x,0)=f(x). \end{equation}

In three dimensions, 
$$ u(x,t)=ct A_tf(x), $$ so, in particular, Stein's spherical maximal theorem implies that 
\begin{equation} \label{waveest} {\left( \int {\left[ \sup_{t>0} \left| \frac{u(x,t)}{t} \right| \right]}^p dx \right)}^{\frac{1}{p}} \leq C_p {||f||}_{L^p({\Bbb R}^3)}, \ \text{for} \ p>\frac{3}{2}, \end{equation} which implies that 
\begin{equation} \label{pointwise} \lim_{t \to 0} \frac{u(x,t)}{t}=f(x) \ \text{if} \ f \in L^p({\Bbb R}^3) \ \text{for} \ p>\frac{3}{2}. \end{equation}  

In higher dimensions similar results are obtained using suitable modifications of the spherical maximal operator. See, for example, \cite{So93} and the references contained therein. 

In order to illustrate how $L^p(\mu)$ spaces, for a suitable measure $\mu$, arise naturally in this context, let us address the following question which was studied in a related context in \cite{A08}. How large is the blowup set of $u(x,1)$ and $\sup_{t \in [1,2]} \frac{u(x,t)}{t}$? More precisely, how large can the Hausdorff dimension of 
\begin{equation} \label{blowupset} \left\{x \in {\Bbb R}^3: u(x,1)=\infty \right\} \ \text{and} \left\{ x \in {\Bbb R}^3: \sup_{t \in [1,2]} \frac{u(x,t)}{t}=\infty \right\} \end{equation} possibly be with $f$ in, say, $L^2({\Bbb R}^3)$? A natural way to approach this problem is to try to prove that 
$$ \int {|u(x,1)|}^2 d\nu(x) \leq C{||f||}^2_{L^2({\Bbb R}^3)}$$ and, correspondingly, that 
$$ \int {\left| \sup_{t \in [1,2]} \frac{u(x,t)}{t} \right|}^2 d\nu(x) \leq C{||f||}^2_{L^2({\Bbb R}^3)}$$ for a measure $\nu$ satisfying the condition $\nu(B(x,r)) \leq Cr^{s_{\nu}}$ with $s_{\nu}>0$ in a suitable range. We could then take $\nu$ to be a Frostman measure on one of the blow-up sets defined in (\ref{blowupset}) above. This problem can (and will) be similarly studied when the initial data $f \in L^p({\Bbb R}^3)$. 

Another context where $L^p(\mu) \to L^p(\nu)$ bounds for classical operators arise is Falconer type problems in geometric measure theory. Falconer proved in 1986 that if the Hausdorff dimension of a compact set $E \subset {\Bbb R}^d$, $d \ge 2$, is $>\frac{d+1}{2}$, then the Lebesgue measure of the distance set $\Delta(E)=\{|x-y|: x,y \in E \}$ is positive. He accomplished this by proving that if $\mu$ is a compactly suported Borel measure satisfying 
$$ I_{\frac{d+1}{2}}(\mu)=\int \int {|x-y|}^{-\frac{d+1}{2}} d\mu(x) d\mu(y)=1$$ then 
\begin{equation} \label{falconer} \mu \times \mu \{(x,y): t \leq |x-y| \leq t+\epsilon \} \leq C(t) \epsilon \end{equation} for any $t>0$. A careful examination of the proof reveals that the Falconer estimate amounts to the classical fact that the spherical averaging operator 
$$ Af(x)=\int_{S^{d-1}} f(x-y) d\sigma(y),$$ where $\sigma$ is the surface measure, maps $L^2({\Bbb R}^d)$ to the Sobolev space $L^2_{\frac{d-1}{2}}({\Bbb R}^d)$. This viewpoint combined with the basic theory of Fourier Integral Operators was used by A. Iosevich, K. Taylor and S. Eswarathasan in \cite{EIT11} to prove an analog of Falconer's result for compact Riemannian manifolds without a boundary. 

Another way to look at (\ref{falconer}) is as an $L^1(\mu) \to L^1(\mu)$ bound for the operator 
\begin{equation} \label{keymamaop} Tf(x)=\sigma*(f\mu)(x) \end{equation} in the special case when $f \equiv 1$. This viewpoint was adopted by K. Taylor in \cite{T12} when she studied the distribution of two-link chains in fractal subsets of ${\Bbb R}^d$. The object analogous to (\ref{falconer}) in this case is 
\begin{equation} \label{taylor} \mu \times \mu \times \mu \{(x,y,z): t \leq |x-z| \leq t+\epsilon; \ t \leq |y-z| \leq t+\epsilon \} \leq C(t)\epsilon^2. \end{equation} 

This estimate can be viewed as the $L^2(\mu) \to L^2(\mu)$ bound for the operator given by (\ref{keymamaop}) in the special case $f \equiv 1$. This viewpoint was further explored in \cite{BIT15} partly based on some of the results obtained in this paper. 

In view of the discussion above, the main focus of this paper is to investigate whether the analog of various classical Lebesgue space bounds hold with the $L^p({\Bbb R}^d) \to L^p({\Bbb R}^d)$ estimate replaced by an $L^p(\mu) \to L^p(\nu)$ bound, where $\mu, \nu$ are compactly supported Borel measures satisfying some natural size assumptions. This leads us to several interesting questions in harmonic analysis, some of which have been explored before, and some that have not. 

The maximal operator estimates in this paper will be studied on a single scale, where the supremum is taken over $t \in [1,2]$. Let us briefly explain why this is essentially unavoidable in the context of general measures. An interesting feature of Stein's spherical maximal operator is that a single-scale bound, where the supremum is taken over $t \in [1,2]$ extends to the multi-scale bound where the supremum is taken over all $t>0$. We shall see that in the context of $L^2(\mu) \to L^2(\nu)$ bounds, the multi-scale maximal estimate is, in general false. We begin with the following simple calculation. Let $\phi$ be a smooth cut-off function supported in the unit ball such that $\widehat{\phi} \ge 0$. Consider the dyadic maximal operator 

$$ {\mathcal D}f(x)=\sup_{j>0} D_jf(x), \ \text{where} \ D_jf(x)=\phi_j*(f\mu)(x), $$ where $\phi_j$ is defined by the relation 
$\widehat{\phi}_j(\xi)=\widehat{\phi}(2^{-j} \xi)$ and $\mu$  is compactly Borel measure satisfying 
\begin{equation} \label{polygrowth} \mu(B(x,r)) \leq Cr^{s_{\mu}}. \end{equation} 

Suppose, in addition, that for every $x$ in the support of $\mu$ and $r$ sufficiently small, 
\begin{equation} \label{lower} \mu(B(x,r)) \ge Cr^{s_{\mu}}. \end{equation} 
Take $f$ such that ${||f||}_{L^2(\mu)}<\infty$ and observe that 
$$ \int D_jf(x) f(x) d\mu(x) \leq {||D_jf||}_{L^2(\mu)} \cdot {||f||}_{L^2(\mu)}.$$ 
On the other hand, the left hand side equals 

\begin{equation} \label{dyadicpizdets} \int {\left|\widehat{f\mu}(\xi)\right|}^2 \widehat{\phi}(2^{-j} \xi) d\xi \ge 
C \int_{|\xi| \leq 2^j} {\left|\widehat{f\mu}(\xi)\right|}^2 d\xi. \end{equation}

A theorem due to Strichartz (\cite{Str90}, Corollary 5.5) says that with $s = s_\mu$ the right hand side of (\ref{dyadicpizdets}) is bounded from below by a constant multiple of $2^{j(d-s)} {||f||}^2_{L^2(\mu)}$. It follows that 
$$ {||{\mathcal D}_jf||}_{L^2(\mu)} \ge C2^{\frac{j(d-s)}{2}}  {||f||}_{L^2(\mu)},$$ so Fatou's lemma implies that 
${|| {\mathcal D}f||}_{L^2(\mu)}=\infty$. Our first result shows that this lower bound extends to the realm of spherical maximal operators but with Lebesgue measure as a target. 
 
\begin{theorem} \label{nehuya} Suppose that $\mu$ is a compactly supported Borel measure that satisfies (\ref{polygrowth}) and (\ref{lower}). Let $f \in L^2(\mu)$ such that ${||f||}_{L^2(\mu)}=1$. Let $\sigma_t$ be defined via its Fourier transform by the formula 
\begin{equation} \label{spheredilate} \widehat{\sigma}_t(\xi)=\widehat{\sigma}(t\xi), \end{equation} where $\sigma$ is the Lebesgue measure on the unit sphere. 

Then 
\begin{equation} \label{dyadicpizdetsest} {|| \sigma_{2^{-j}}*(f\mu)||}_{L^2({\Bbb R}^d)} \ge C2^{\frac{j(d-s)}{2}} \end{equation} and, consequently, 

\begin{equation} \label{divergencepizdets} {\left|\left| \sup_{j>0} | \sigma_{2^{-j}}*(f\mu)| \right|\right|}_{L^2({\Bbb R}^d)}=\infty. \end{equation} 

\end{theorem} 

To see this, observe that by Plancherel, 
\begin{equation} \label{l2r} {|| \sigma_{2^{-j}}*(f\mu) ||}^2_{L^2({\Bbb R}^d)}=
\int {|\widehat{f\mu}(\xi)|}^2 {|\widehat{\sigma}(2^{-j} \xi)|}^2 d\xi. \end{equation} 

If $|\xi| \leq c2^j$ with $c$ small enough, then ${|\widehat{\sigma}(2^{-j} \xi)|}^2 \ge c'>0$. It follows that the expression above is bounded from below by 

$$ \int_{|\xi| \leq c2^j} {\left|\widehat{f\mu}(\xi)\right|}^2 d\xi.$$ 

Invoking Theorem 5.5 in \cite{Str90} once again, we see that this expression is bounded from below by 
$$ C 2^{j(d-s)} {||f||}_{L^2(\mu)}^2,$$ which shows that the operator norm of 
$\sigma_{2^{-j}}*(f\mu)$ from $L^2(\mu) \to L^2({\Bbb R}^d)$ is at least $C 2^{\frac{j(d-s)}{2}}$. This completes the proof of Theorem \ref{nehuya}. 

\qed

\vskip.125in 

\begin{remark} There is no hope of establishing a "universal" negative result like Theorem \ref{nehuya} when $s_{\nu}<d$. In order to see this consider the case $d \ge 3$, $s_{\mu}=s_{\nu}=s$ and $\mu=\nu$ is the surface measure on the unit sphere. In this case, the bound 
$$ {\left( \int \sup_{t>0} {|\sigma_t*(f\mu)(x)|}^p d\mu(x) \right)}^{\frac{1}{p}} \leq C_p {\left( \int {|f(x)|}^p d\mu(x) \right)}^{\frac{1}{p}} \ \text{for} \ p>\frac{d-1}{d-2}$$ holds as a consequence of the known extension of the Bourgain-Stein result to compact Riemannian manifolds without a boundary. See \cite{So93}, Chapter 7, and the references contained therein. Of course the difference between \eqref{divergencepizdets} and the inequality above is that the target norm above is weighted.


\end{remark} 

\vskip.125in

Having established the necessity of working away from the small times scales in (\ref{divergencepizdets}) above, we now turn our attention to positive results. This paper is organized as follows. \begin{itemize} 

\item Before turning to maximal operators, we investigate $L^2(\mu) \to L^2(\nu)$ bounds for the convolution operators 
\begin{equation}  \label{convfractal} T_{\lambda^{\epsilon}}(f\mu)(x)=\lambda^{\epsilon}*(f\mu)(x), \end{equation} where $\lambda$ is a tempered distribution on ${\Bbb R}^d$ and $\lambda^{\epsilon}=\lambda*\rho_{\epsilon}(x)$, with $\rho_{\epsilon}(x)=\epsilon^{-d} 
\rho(x/ \epsilon)$, $\rho \in C^{\infty}_0({\Bbb R}^d)$ and $\int \rho=1$. This analysis is carried out in Section \ref{convolutionfractal} below. We will also see that our results easily extend to a wide class of operators satisfying $L^2$-Sobolev bounds and a natural micro-local assumption.  

\item We shall then prove $L^p(\mu) \to L^p(\nu)$ bounds for the maximal operator 
\begin{equation} \label{fractalsphericalmax} {\mathcal A}_{\mu}f(x)=\sup_{t \in [1,2]} \sigma_t*(f\mu)(x), \end{equation} where $\sigma_t$ is defined in (\ref{spheredilate}) above. This analysis is carried out in Section \ref{fiomaximal} below. As the reader shall see, the estimates rely on robust Fourier Integral Operator estimates and thus the results extend far beyond spherical maximal operators. 

\item In the next section, we shall describe some applications of the results we have obtained to the initial problem for the wave equation in three dimensions. In this context we also obtained detailed $L^p(\mu) \to L^p(\nu)$ estimates for the spherical averaging operator $A_tf(x)=\int_{S^{d-1}} f(x-y) d\sigma(y)$. 

\item Throughout the paper we emphasize the following consequence of our $L^p(\mu) \to L^p(\nu)$ bounds. By setting $\mu$ to be the Lebesgue measure and $\nu$ to be a Frostman measure on the set $E_f=\{x: Tf(x)=\infty\}$, we obtain upper bounds for the Hausdorff dimension of blowup sets for the operators under consideration. Connections with the local smoothing conjecture are also explored. 

\end{itemize}

\vskip.25in 

\section{Classical convolution inequalities in a fractal setting} 
\label{convolutionfractal} 

\vskip.125in 

In this section we show that bounds for these and other convolution operators extend to the setting where $L^p({\Bbb R}^d)$ is replaced by $L^p(\mu)$, where $\mu$ is a compactly supported Borel measure satisfying 
\begin{equation} \label{upperad} \mu(B(x,r)) \leq Cr^{s_{\mu}} \end{equation} for some $d\geq s_{\mu}>0$ and every $x \in supp(\mu)$ and $r \in [0, diam(supp(\mu))]$. 

Let $\lambda$ be a tempered distribution, and denote by $\lambda^{\epsilon}$ the convolution of $\lambda$ with $\rho_{\epsilon}(x) \equiv \epsilon^{-d} \rho \left( \frac{x}{\epsilon} \right)$, $\epsilon>0$, with $\rho \in C_0^{\infty}({\Bbb R}^d)$, and $\int \rho(x) dx=1$. Then $\lambda^{\epsilon}$ is a $C^{\infty}({\Bbb R}^d)$ function. Define 

\begin{equation} \label{operatordef} T_{\lambda^{\epsilon}}f(x)=\int \lambda^{\epsilon}(x-y) f(y) d\mu(y), \end{equation} where $\mu$ is a compactly supported Borel measure satisfying (\ref{upperad}) above, and $\lambda$ is a tempered distribution whose Fourier transform is a locally integrable function satisfying 
\begin{equation} \label{decay} |\widehat{\lambda}(\xi)| \leq C{|\xi|}^{-\alpha} \end{equation} for some $\alpha \in \left[0, \frac{d}{2} \right)$. 

If $\mu$ is the Lebesgue measure on ${\Bbb R}^d$, the Plancherel theorem says that the $L^2({\Bbb R}^d)$ bound of $T_{\lambda^{\epsilon}}$ holds if and only if $\widehat{\lambda}$ is bounded. If $\mu$ is not the Lebesgue measure, Plancherel is not available. As a substitute, we have the following result. 

\begin{theorem} \label{main} Let $T_{\lambda^{\epsilon}}$ be as in (\ref{operatordef}) above with $\lambda$ satisfying (\ref{decay}). Let $\mu$ be a compactly supported Borel measure satisfying (\ref{upperad}) and suppose that $\nu$ is another compactly supported Borel measure satisfying (\ref{upperad}) with $s_{\mu}$ replaced by $s_{\nu}$. 

\vskip.125in 

\begin{itemize} \item i) Let $s=\frac{s_{\mu}+s_{\nu}}{2}$. Suppose that $\alpha>d-s$. Then 
\begin{equation} \label{ytm} {||T_{\lambda^{\epsilon}}f||}_{L^2(\nu)}  \leq C{||f||}_{L^2(\mu)} \end{equation} with constant $C$ independent of $\epsilon$. 

\vskip.125in

\item ii) If $\alpha \leq d-s$, then $T^{\epsilon}_{\lambda}$ does not in general map $L^2(\mu)$ to $L^2(\nu)$ with constants independent of $\epsilon>0$. 

\end{itemize} 

\end{theorem} 

\begin{remark} \label{noepsilon} Note that if $\lambda$ is, say, a compactly supported finite Borel measure, then the $\epsilon$-mollification above is not necessary. Moreover, we can interpret (\ref{ytm}) without $\epsilon$ in the standard way. The tempered distribution $\lambda*(f\mu)$ is a weak limit of the $C^{\infty}$ functions $\lambda^{\epsilon}*(f\mu)$ that are uniformly bounded in $L^2(\nu)$ by (\ref{ytm}). Thus the tempered distribution $T_{\lambda}f \equiv \lambda*(f \mu)$ can be interpreted as an $L^2(\nu)$ function and the $\epsilon$ can be removed.  \end{remark} 

\vskip.125in

\begin{remark} \label{cooolops} We stated Theorem \ref{main} for convolution operators for the sake of convenience and notational simplicity. An interested reader can easily check that the same conclusion holds for the operators satisfying the following assumptions: \begin{itemize} 

\item i) We have 
$$ Tf(x)=\int K(x,y) f(y) d\mu(y),$$ where $K$ is a measurable function on ${\Bbb R}^d \times {\Bbb R}^d$. 

\item ii) The operator 
$$ T^{E}f(x)=\int_{{\Bbb R}^d} K(x,y) f(y) dy \ \text{boundedly maps} \ L^2({\Bbb R}^d) \to L^2_{\alpha}({\Bbb R}^d).$$ 

\item iii) There exists a finite $C>0$ such that 
$$ support(\widehat{Tf}) \subset C \cdot support(\widehat{f})=\{Cx: x \in support(\widehat{f})\}.$$ 

\end{itemize} 

We note that these assumptions are, in particular, satisfied by non-degenerate Fourier Integral Operators of order $-\alpha$. See, for example, \cite{So93} and the references contained therein. 

\end{remark}

\begin{corollary}\label{BlowUpSetforThmMain}
 Let $T_{\lambda}$  and $\mu$ be as in Theorem \ref{main} and Remark \ref{noepsilon}. Then the function $T_{\lambda}f \equiv \lambda*(f \mu)$ is actually defined up to a set of Hausdorff $s_{\nu}$-measure zero, and moreover, the set $\left\{   x \in \mathbb{R}^d :  T_{\lambda}f (x) \equiv \lambda*(f \mu) (x) = \infty   \right\}$ also has Hausdorff $s_{\nu}$-measure zero. More generally, the same result holds for operators described in Remark \ref{cooolops}. 
\end{corollary}

We shall discuss the sharpness of this result in the concrete case when 
$$T_{\lambda}f(x)=Af(x)=\int_{S^{d-1}} f(x-y) d\sigma(y),$$ where $\sigma$ is the surface measure on $S^{d-1}$ in the last section of the paper. 

\vskip.125in 

\subsection{Proof of Corollary \ref{BlowUpSetforThmMain}}

Let the set of non-definition (or the blow-up set) be denoted by $E$ and denote the Hausdorff $ s_{\nu}$-measure by $\mathcal{H}^{s_{\nu}}$. Then if $\mathcal{H}^{s_{\nu}} (E) > 0$,  $E$ would contain a compact subset $F$ also with $\mathcal{H}^{s_{\nu}} (F) > 0$. By Frostman's lemma, there would exist a Borel measure $\nu$ supported on $F$, satisfying  (\ref{upperad}) with $s_{\mu}$ replaced by $s_{\nu}$. That measure $\nu$ would provide a contradiction to the Theorem \ref{main}.

\qed

\subsection{Proof of Theorem \ref{main} , part i)} It is enough to show that if $g \in L^2(\nu)$, then 
\begin{equation} \label{setup} |<T_{\lambda^{\epsilon}}f, g\nu>| \ \leq C{||f||}_{L^2(\mu)} \cdot {||g||}_{L^2(\nu)}, \end{equation} where the constant $C$ is independent of $\epsilon$. 

\vskip.125in 

The left hand side of (\ref{setup}) equals 
\begin{equation} \label{preinversion} \int \lambda^{\epsilon}*(f\mu)(x) g(x) d\nu(x). \end{equation} 

Indeed, 
$$ \lambda^{\epsilon}*(f\mu)(x)=\int e^{2 \pi i x \cdot \xi} \widehat{\lambda}(\xi) \widehat{\rho}(\epsilon \xi)
\widehat{f\mu}(\xi) d\xi$$ for every $x \in {\Bbb R}^d$ because the left hand side is a continuous $L^2({\Bbb R}^d)$ function and $$\widehat{\lambda}(\cdot) \widehat{\rho}(\epsilon \cdot) \widehat{f\mu}(\cdot) \in L^1 \cap L^2({\Bbb R}^d).$$

It follows that (\ref{preinversion}) equals 

$$\int \int e^{2 \pi i x \cdot \xi} \widehat{\lambda}(\xi) \widehat{\rho}(\epsilon \xi) \widehat{f \mu}(\xi) d\xi g(x) d\nu(x).$$ Applying Fubini, we see that this expression equals 
$$ \int \int e^{2 \pi i x \cdot \xi} g(x) d\nu(x) \widehat{\lambda}(\xi) \widehat{\rho}(\epsilon \xi) \widehat{f \mu}(\xi) d\xi$$
$$=\int \widehat{\lambda}(\xi) \widehat{\rho}(\epsilon \xi) \widehat{f \mu}(\xi)  \widehat{g \nu}(\xi) d\xi.$$

The modulus of this expression is bounded by an $\epsilon$-independent constant multiple of 

$$ \int {|\xi|}^{-\alpha} |\widehat{f\mu}(\xi)| \cdot  |\widehat{g \nu}(\xi)| d\xi.$$ 

By Cauchy-Schwartz, this expression is bounded by 

\begin{equation} \label{almostthere} {\left( \int {|\widehat{f\mu}(\xi)|}^2 {|\xi|}^{-\alpha_{\mu}} d\xi \right)}^{\frac{1}{2}} \cdot {\left( \int {|\widehat{g\nu}(\xi)|}^2 {|\xi|}^{-\alpha_{\nu}} d\xi \right)}^{\frac{1}{2}}=\sqrt{A} \cdot \sqrt{B}, \end{equation} where $\alpha_{\mu}, \alpha_{\nu}>0$ and $\frac{\alpha_{\mu}+\alpha_{\mu}}{2}=\alpha$. 

\vskip.125in 

\begin{lemma} \label{Aestimation} With the notation above, we have 

$$ A \leq C{||f||}^2_{L^2(\mu)}; B \leq C{||g||}^2_{L^2(\nu)}$$ if 
\begin{equation} \label{numuconditions} \alpha_{\nu}>d-s_{\nu}, \alpha_{\mu}>d-s_{\mu}, \ \text{respectively}. \end{equation} 
\end{lemma} 

\vskip.125in 

Lemma \ref{Aestimation} can be deduced, via a dyadic decomposition, from the following fact due to Strichartz (\cite{Str90}). With the notation above, 
\begin{equation} \label{bob} \left| \sup_{r \ge 1} r^{-(d-s_{\mu})} \int_{B(x,r)} {|\widehat{f\mu}(\xi)|}^2 d\xi \right| \leq C{||f||}^2_{L^2(\mu)}. \end{equation} 

Instead, we give a direct argument in the style of the proof of Theorem 7.4 in \cite{W04}. It is enough to prove the estimate for $A$ since the estimate for $B$ follows from the same argument. By Proposition 8.5 in \cite{W04}, 
\begin{equation} \label{schursetup} A=\int \int f(x)f(y) {|x-y|}^{-d+\alpha_{\mu}} d\mu(x)d\mu(y)=<f,Uf>_{L^2(\mu)}, \end{equation} where 

$$ Uf(x)=\int {|x-y|}^{-d+\alpha_{\mu}} f(y) d\mu(y).$$ 

Observe that (we assume for simplicity that the diameter of the support of $\mu$ is $\leq 1$):
$$ \int {|x-y|}^{-d+\alpha_{\mu}} d\mu(y)= \int {|x-y|}^{-d+\alpha_{\mu}} d\mu(x)$$
$$ \leq C \sum_{j>0} 2^{j(d-\alpha_{\mu})} \int_{2^{-j} \leq |x-y| \leq 2^{-j+1}} d\mu(y)$$
$$ \leq C' \sum_{j>0} 2^{j(d-\alpha_{\mu}-s_{\mu})} \leq C'' \ \text{if} \ \alpha_{\mu}>d-s_{\mu}.$$

It follows by Schur's test (see Lemma 7.5 in \cite{W04} and the original argument in \cite{Schur11}) that 
$${||Uf||}_{L^2(\mu)} \leq C'' {||f||}_{L^2(\mu)}$$ and we are done in view of (\ref{schursetup}) and the Cauchy-Schwartz inequality. 

\vskip.25in 

\subsection{Proof of  Theorem \ref{main}, part ii)} 

We shall consider the case $s_{\mu}=s_{\nu}=s$, but an interested reader can easily generalize this example. Let $\lambda(x)={|x|}^{-d+\alpha}\chi_B(x)$, where $B$ is the unit ball, and suppose that $\mu$ is the restriction of the $s$-dimensional Hausdorff measure to a compact Ahlfors-David regular set of dimension $s$. Then 

$$ Tf(x)=\int {|x-y|}^{-d+\alpha} f(y) d\mu(y).$$ 

Let $f \equiv 1$ and observe that 

$$ T1(x) \approx \sum_{j\geq -1} 2^{j(d-\alpha)} \int_{2^{-j} \leq |x-y| \leq 2^{-j+1}} d\mu(y)$$

$$ \gtrsim \sum_{j\geq -1} 2^{j(d-\alpha)} \cdot \mu(B(x,2^{-j})) \approx \sum_{j\geq -1} 2^{j(d-(s+\alpha))}$$ and this quantity is infinite if $s \leq d-\alpha$. 

\qed

\vskip.25in 

\section{Spherical maximal operator over fractals} 
\label{fiomaximal}

\vskip.125in 

As we noted above, the classical Bourgain-Stein maximal operator is bounded on $L^p({\Bbb R}^d)$, $d \ge 2$, for $p>\frac{d}{d-1}$. An unusual feature of the fractal analog of this result, which we are about to state, is that the range of $p$s for which the maximal operator is bounded when the supremum is taken over a "single scale" $t \in [1,2]$ is a bounded interval. This is a reflection of the fact that if $\mu$ satisfies (\ref{polygrowth}) with $s_{\mu}<d$, then $\sigma*\mu$ is not, in general, bounded. 

\begin{theorem} \label{maximalmain} Let $\mu, \nu$ denote compactly supported Borel measures satisfying the condition (\ref{polygrowth}). Let $A_tf(x)=\sigma_t*(f\mu)(x)$, where $\sigma_t$ is the spherical measure defined as above. Let 
${\mathcal A}f(x)=\sup_{t \in [1,2]} |A_tf(x)|$. 

\vskip.125in 

i) Suppose that $s_{\mu}+s_{\nu}>d+2$ and $s_{\mu}>1$. Then 
\begin{equation} \label{bound1} {\mathcal A}: L^p(\mu) \to L^p(\nu)  \ \text{for} \ p \in \left(\frac{d+s_{\mu}-s_{\nu}}{s_{\mu}-1}, p_{U} \right), \end{equation} where $p_{U}$ is the supremum of the set of $p$s such that 
\begin{equation} \label{upperrange} d-s_{\mu}+\frac{s_{\mu}-s_{\nu}}{p}<\frac{d-2}{p}+\epsilon(p), \end{equation} where 
$$ \epsilon(p)=\frac{1}{2p} \ \text{if} \ p \ge 4 \ \text{and} \ \epsilon(p)=\frac{1}{2}\left( \frac{1}{2}-\frac{1}{p} \right) \ \text{if} \ 2 \leq p \leq 4.$$ 

\vskip.125in 

ii) Suppose that $s_{\mu}+s_{\nu} = d+2$, $s_{\mu}>1$ and $2 \leq s_{\nu}< \frac{9}{4}$. Then 
\begin{equation} \label{bound3} {\mathcal A}: L^p(\mu) \to L^p(\nu)  \ \text{for} \ p \in (2, 4]. \end{equation} 
\vskip.125in 

iii) Suppose that $s_{\mu}+s_{\nu} < d+2$, and $ 3s_{\mu}+s_{\nu} > 3d + \frac{3}{2}$. Then 
\begin{equation} \label{bound4} {\mathcal A}: L^p(\mu) \to L^p(\nu)  \ \text{for} \ p \in \left( \frac{ s_{\nu} - s_{\mu} +d - \frac{5}{2} }{ d - s_{\mu} - \frac{1}{4} } ,  \frac{ s_{\nu} - s_{\mu} +d - \frac{3}{2} }{ d - s_{\mu} } \right). \end{equation} Note that in this case $2$ is not in the interval, but $4$ is in the interval.

\end{theorem} 

\vskip.125in 

\begin{remark} Observe that if $s_{\nu}=d$, the left endpoint of the interval in i) above is $\frac{s_{\mu}}{s_{\mu}-1}$. We shall see in Theorem \ref{steinexample} below that this estimate is best possible. The right endpoint is probably not sharp since the proof relies on the known local smoothing estimates which are not believed to be sharp. Recall (see e.g. \cite{So93}) that the sharp local smoothing results in this context would imply the celebrated Kakeya conjecture. \end{remark} 

\vskip.125in 

\begin{remark} 
It is particularly interesting to illustrate part iii) of Theorem \ref{maximalmain} in the case $d=2$ where the parts i) - ii) condition $s_{\mu}+s_{\nu}\geq d+2$ is never satisfied. Let $s_{\mu}=s_{\nu}=s$. The condition $ 3s_{\mu}+s_{\nu} > 3d + \frac{3}{2}$ reduces to $s > \frac{15}{8}$, and the  $L^p$ range of bounds is  $\left( \frac{1}{2 \left(s-\frac{7}{4} \right)} ,\frac{1}{2(2-s)} \right)$.
\end{remark} 

\vskip.125in 

\begin{remark} 
The other cases not covered in the statement of Theorem \ref{maximalmain} yield contradictions with the conditions of convergence in the method of proof of the theorem.
\end{remark} 

\vskip.125in 

The following corollary of Theorem \ref{maximalmain} is the corresponding analogue of Corollary \ref{BlowUpSetforThmMain}, and is proved in the same way. We state it in 3 parts, each corresponding to the part with the same numbering in Theorem \ref{maximalmain}. Note that, as a consequence of it, one gets control on the blowup set (both in $x$ for a fixed $t$ and jointly in $(x,t)$) of solutions to PDEs whenever the solution is controlled by the maximal operator ${\mathcal A}f(x)$, e.g. as in the case of the wave equation, as we will see soon. More explicitly, the following corollary gives, for a given datum $f \in L^p(\mu)$, control on the size of the blowup set of the solution to the wave equation for that particular datum $f$. This of course includes the case when $\mu$ is Lebesgue measure in $\mathbb{R}^d$, in which case $s_\mu = d$.

\vskip.125in

\begin{corollary}\label{BlowUpSetforThmMaximalMain}
 Let $\mu$ denote a compactly supported Borel measure satisfying the condition (\ref{polygrowth}).  Let ${\mathcal A}f(x)$ be as in Theorem \ref{maximalmain}, for a function $f \in L^p(\mu)$, for some $p > 1$. Let $p_f$ be the supremum of the set of $p$'s so that $f \in L^p(\mu)$. 

\vskip.125in 

i) Assume $s_\mu > 1$, and that  $p_f > 2$. Let $s_1$ be defined by
\begin{equation} \label{upperrangeForCorollary} 
d-s_{\mu}+\frac{s_{\mu}-s_1}{p_f} = \frac{d-2}{p_f}+\epsilon(p_f), 
\end{equation} 
where $\epsilon(p)$ is as in \eqref{upperrange}, and let $s_f = \max\left\{s_1, d+2-s_\mu \right\}$. 
Then the blowup set 
$$ E_f := \left\{ x \in \mathbb{R}^d : {\mathcal A}f(x) = \infty   \right\} $$ 
satisfies $\dim_{\mathcal H} (E_f) \leq s_f$, where $\dim_{\mathcal H}$ denotes Hausdorff dimension.

Moreover, if $s_f > d+2-s_\mu$, then $\mathcal{H}^{s_f} (E_f) = 0$, where $\mathcal{H}^{s_f}$ denotes Hausdorff $s_f$-dimensional measure.

\vskip.125in 

ii) Assume $d - \frac{1}{4} < s_\mu \leq d$, and that  $p_f > 2$. Then $\mathcal{H}^{d+2-s_\mu} (E_f) = 0$.

\vskip.125in 

iii) Assume $d - \frac{1}{4} < s_\mu \leq d$, and that  $p_f > 4$. Then  $\dim_{\mathcal H} (E_f) \leq 3(d - s_\mu) + \frac{3}{2}$. Note that if $d - \frac{1}{4} < s_\mu$, then $ 3(d - s_\mu) + \frac{3}{2} < d+2-s_\mu$.

\end{corollary}

\vskip.125in 

\subsection{Proof of Corollary \ref{BlowUpSetforThmMaximalMain}}

We will use repeatedly without mentioning it that since $\mu$ is compactly supported, $f \in L^p(\mu)$ for all $p < p_f$.

\vskip.125in 

i) If $s_1 > d+2-s_\mu$, then $s_f = s_1$. Then, in \eqref{bound1}, we see that $\frac{d+s_{\mu}-s_f}{s_{\mu}-1} < \frac{2s_{\mu}-2}{s_{\mu}-1} = 2$, and by \eqref{upperrangeForCorollary}, $p_{U} \geq p_f > 2$. As a consequence, for any compactly supported measure $\nu$ satisfying (\ref{upperad}) with $s_{\mu}$ replaced by $s_{\nu}$, we get from \eqref{bound1} that ${\mathcal A}: L^2(\mu) \to L^2(\nu) $. Then if $\mathcal{H}^{s_f} (E_f) > 0$, by Frostman's lemma, there exists $\nu$ compactly supported in $E_f$ satisfying (\ref{upperad}) with $s_{\mu}$ replaced by $s_{\nu}$. This $\nu$ provides a contradiction with Theorem \ref{maximalmain} i). Consequently, $\mathcal{H}^{s_f} (E_f) = 0$.

If $s_1 = d+2-s_\mu$, then $s_f = s_1$.  Then, in \eqref{bound1}, we see that $\frac{d+s_{\mu}-s_f}{s_{\mu}-1} = \frac{2s_{\mu}-2}{s_{\mu}-1} = 2$. Consider $\widetilde{s_f} = s_f + \delta$, for sufficiently small $\delta > 0$. Now  $s_\mu + \widetilde{s_f} > d+2$, and  $\frac{d+s_{\mu}-\widetilde{s_f}}{s_{\mu}-1} < \frac{d+s_{\mu}-s_f}{s_{\mu}-1} = 2$. Also from \eqref{upperrangeForCorollary} we see that $d-s_{\mu}+\frac{s_{\mu}-\widetilde{s_f}}{p_f} < \frac{d-2}{p_f}+\epsilon(p_f)$, thus, for $\widetilde{s_f}$, we get that $p_{U} > p_f > 2$. Then as in the previous case, we get $\mathcal{H}^{\widetilde{s_f}} (E_f) = 0$. Sending $\delta \to 0$, we get the desired conclusion.

If $s_1 < d+2-s_\mu$, then $s_f = d+2-s_\mu$. Exactly the same reasoning as the previous case (when  $s_1 = d+2-s_\mu$) applies to give the desired conclusion.

\vskip.125in 

ii) Take $s_v = d + 2 - s_\mu$, under the stated conditions we get that $2 \leq s_{\nu}< \frac{9}{4}$, and apply Theorem \ref{maximalmain}, part ii) together with the same reasoning as part i).

\vskip.125in 

iii) First note that  $ 3(d - s_\mu) + \frac{3}{2} < d+2-s_\mu$ if and only if $d - \frac{1}{4} < s_\mu$. Now choose $s_\nu = 3(d - s_\mu) + \frac{3}{2} + \delta$ for sufficiently small $\delta >0$, and apply Theorem \ref{maximalmain}, part iii) together with the same reasoning as part i), to conclude that  $\dim_{\mathcal H} (E_f) \leq s_\nu$. Now send $\delta \to 0$.

\qed

\vskip.125in 

We now address the sharpness of Theorem \ref{maximalmain}, at least to a certain extent. 

\vskip.125in 

\begin{theorem} \label{steinexample} Let ${\mathcal A}$ be defined as above. Then ${\mathcal A}: L^p(\mu) \to L^p(\nu)$ does not in general hold if $s_{\mu}>1$ and $1 < p \leq \frac{s_{\mu}}{s_{\mu}-1}$, or if $0 \leq s_{\mu} \leq1$ and $1 < p < \infty$.
\end{theorem} 

\begin{theorem} \label{mattilaexample} Let $s_{\mu}, \mu, \nu$ be defined as above. 

\vskip.125in 

a) Suppose that $d \ge 2$ and $s_{\mu}<1+\frac{2}{p}$, with $p \in [1,2]$. Then ${\mathcal A}: L^p(\mu) \to L^p(\nu)$ does not in general hold. Note that this exponent $p$ matches the left endpoint in part i) of Theorem \ref{maximalmain} if $d=2$ and $s_{\mu} = s_{\nu}$. 

b) Suppose that $p \ge 2$, $d=2$ and $s_{\mu}<\frac{3-\frac{2}{p}}{2-\frac{2}{p}}$. Then ${\mathcal A}: L^p(\mu) \to L^p(\nu)$ does not in general hold. In particular, if $p \ge 2$, ${\mathcal A}: L^p(\mu) \to L^p(\nu)$ does not in general hold for any $p \ge 2$ if $s_{\mu}<\frac{3}{2}$. 

c) Suppose that $d \ge 3$ and $p \ge 2$. Then ${\mathcal A}: L^p(\mu) \to L^p(\nu)$ does not in general hold for any $p \ge 2$ if $s_{\mu}<2$. 

\end{theorem} 

\vskip.125in 

\vskip.25in 

\subsection{Proof of Theorem \ref{maximalmain}} 

\vskip.125in 

We shall make use of the following fundamental results, due to Seeger, Sogge and Stein (\cite{SSS91}) and Mockenhaupt, Seeger, Sogge (\cite{MSS92}), respectively. See \cite{So93}, Theorem 6.2.1 and Theorem 7.1.1 for the description of the notation and the background. 

\begin{theorem} Let $X,Y$ be $d$-dimensional $C^{\infty}$ manifolds and let ${\mathcal F} \in I^m(X,Y, {\mathcal C})$, with ${\mathcal C}$ being a locally the graph of a canonical transformation. Then 
$$ {\mathcal F}: L^p_{comp}(Y) \to L^p_{loc}(X) \ \text{if} \ 1<p<\infty \ \text{and} \ m \leq -(d-1) \left| \frac{1}{p}-\frac{1}{2} \right|.$$ 
\end{theorem} 

\vskip.125in 

\begin{theorem} Suppose that ${\mathcal F} \in I^{m-\frac{1}{4}}(Z,Y; {\mathcal C})$ where ${\mathcal C}$ satisfies the non-degeneracy assumptions and the cone condition above. Then ${\mathcal F}: L^p_{comp}(Y) \to L^p_{loc}(Z)$ if $m \leq -(d-1) \left| \frac{1}{p}-\frac{1}{2} \right|+\epsilon(p)$, where $\epsilon(p)=\frac{1}{2p}$ if $p \ge 4$ and $\epsilon(p)=\frac{1}{2} \left( \frac{1}{2}-\frac{1}{p} \right)$ if 
$2 \leq p \leq 4$. 
\end{theorem} 

\vskip.125in 

%
%
%
%

By the Sobolev embedding theorem, it is enough to show that for any $\delta>0$
\begin{equation} \label{set} 
{\left( \int \int_1^2 {\left|{\left(\frac{d}{dt}\right)}^{\frac{1}{p}+\delta} A_tf(x)\right|}^p dt d\nu(x) \right)}^{\frac{1}{p}} 
\leq C{||f||}_{L^p(\mu)}. \end{equation} 

To prove this, it is enough to show that if $g \in L^{p'}(\nu \times \gamma)$, where $\gamma$ is the Lebesgue measure on $[1,2]$, then 
\begin{equation} \label{polarkaif} \int_1^2 \int {\left(\frac{d}{dt}\right)}^{\frac{1}{p}+\delta}A_tf(x) g(x,t) d\nu(x) dt \leq C {||f||}_{L^p(\mu)} \cdot {||g||}_{L^{p'}(\nu \times \gamma)}. \end{equation} 

Define the Littlewood-Paley projection by the relation 
$$ \widehat{P_jf}(\xi)=\widehat{f}(\xi) \beta(2^{-j} \xi),$$ where $\beta$ is a smooth cut-off function supported in the annulus 
$\{\xi: 1/2 \leq |\xi| \leq 4 \}$, identically equal to $1$ in the annulus $\{\xi: 1 \leq |\xi| \leq 2 \}$ and $\{\sum_{j\geq1} \beta(2^{-j} \xi) \} + \beta_0(\xi) \equiv 1$, where $\beta_0$ is another smooth cutoff function supported in $\{\xi: |\xi| \leq 2 \}$. $P_0$ will denote the corresponding Littlewood-Paley projection for $\beta_0$.

Define 
$$ A_t^jf(x)=\sigma_t*(P_j(f\mu))(x).$$ 

The left hand side of (\ref{polarkaif}) equals 

$$ \sum_j \int_1^2 \int {\left(\frac{d}{dt}\right)}^{\frac{1}{p}+\delta}A_t^jf(x) g(x,t) d\nu(x) dt.$$ 

The case $j = 0$ is easily handled using the Hardy-Littlewood maximal function, so we confine our attention to the positive frequencies. Observe that the sum above essentially equals 
\begin{equation} \label{lplocalized} \sum_j \int_1^2 \int {\left(\frac{d}{dt}\right)}^{\frac{1}{p}+\delta}A_t^jf(x) P_j(g_t\nu)(x) dx dt,\end{equation} where $g_t(x)=g(x,t)$. 

More precisely, we must consider $P_k(f\mu)$ with $|k-j| \leq 2$, but this reduces to the expression above by the triangle inequality and relabeling. 

Applying H\"older to (\ref{lplocalized}) yields 
\begin{equation} \label{platzdarm} \sum_j {\left( \int_1^2 \int {\left| {\left(\frac{d}{dt}\right)}^{\frac{1}{p}+\delta}A_t^jf(x) \right|}^p dx dt \right)}^{\frac{1}{p}} \cdot {\left( 
\int_1^2 \int {\left| P_j(g_t \nu)(x) \right|}^{p'} dx dt \right)}^{\frac{1}{p'}}=\sum_j I \cdot II.  \end{equation} 


%

Lemma \ref{Aestimation} gives us a bound for the term $II$ in (\ref{platzdarm}) in the case $p=2$. We obtain a bound for other exponents by interpolation. We have 
$$ P_j(g_t\nu)(x)=2^{dj} \widehat{\beta}(2^j \cdot)*(g_t\nu)(x)$$
\begin{equation} \label{Linfty}=2^{dj} \int \widehat{\beta}(2^j(x-y)) g_t(y) d\nu(y) \leq C2^{j(d-s_{\nu})} {||g_t||}_{L^{\infty}(\nu)}. \end{equation}

Similarly, 

\begin{equation} \label{L1} \int P_j(g_t\nu)(x) dx \leq C2^{dj} \int \int |\widehat{\beta}(2^j(x-y))| |g_t(y)| d\nu(y) dx \leq 
C' {||g_t||}_{L^1(\nu)}. \end{equation} 

It follows that if $1 \leq p \leq 2$, for any $\eta >0$,
\begin{equation} \label{Lpprime} {\left(\int_1^2 \int {|P_j(g_t\nu)(x)|}^{p'} dx dt \right)}^{\frac{1}{p'}} \leq C_\eta 2^{j \frac{d-s_{\nu} + \eta}{p}} {||g||}_{L^{p'}(\nu \times \gamma)} \end{equation} by interpolating between (\ref{Linfty}) and the following variant of Lemma \ref{Aestimation}. 

\begin{lemma} \label{Aestimationvar} With the notation above, for any $\eta>0$,
$$ \left| \int {|\widehat{f \nu}(\xi)|}^2 \beta(2^{-j}\xi) d\xi \right| \leq C_\eta 2^{j(d-s_{\nu} + \eta)} {||f||}_{L^2(\nu)}^2.$$ 
\end{lemma} 

The proof of Lemma \ref{Aestimationvar} is the same as the proof of Lemma \ref{Aestimation} (and actually it can be easily deduced from it). 

\vskip.125in 

If $2 \leq p \leq \infty$, (\ref{Lpprime}) holds by interpolating between Lemma \ref{Aestimationvar} and (\ref{L1}). This completes the estimation of the term $II$ in (\ref{platzdarm}). 

\subsubsection{Proof of part i)} We now turn our attention to the estimation of the term $I$ in the regime $s_{\mu}+s_{\nu}>d+2$. If $1 \leq p \leq 2$, we use the Seeger-Sogge-Stein bound (\cite{SSS91}). This yields 
\begin{equation} \label{estI} I \leq C 2^{-j \frac{d-1}{2}} \cdot 2^{j (d-1)\left| \frac{1}{2}-\frac{1}{p}\right|} 2^{\frac{j}{p}} 2^{j \delta}
{||P_j(f\mu)||}_{L^p({\Bbb R}^d)}, \end{equation} where the factor $2^{\frac{j}{p}}2^{j \delta}$ comes from the differentiation in 
$t$ of order $\frac{1}{p} + \delta$. 

Using (\ref{Lpprime}) this, in turn, is 
$$ \leq C_\eta 2^{-j \frac{d-1}{p'}} 2^{\frac{j}{p}} 2^{j \delta} 2^{j \frac{d-s_{\mu} + \eta}{p'}}{||f||}_{L^p(\mu)}.$$

Plugging this into (\ref{platzdarm}) along with (\ref{Lpprime}), we obtain 

$$ I \cdot II \leq C_\eta 2^{-j \frac{d-1}{p'}} 2^{\frac{j}{p}} 2^{j (\delta + \eta)} 2^{j \frac{d-s_{\mu}}{p'}} 2^{j \frac{d-s_{\nu}}{p}}{||f||}_{L^p(\mu)} {||g||}_{L^{p'}(\nu \times \gamma)}.$$ 

Summing in $j$ and recalling that $\delta$ and $\eta$ are arbitrarily small, we see that the geometric series converges if 
\begin{equation}\label{ConditionForConvergenceIfPLessThan2} 
p>\frac{d+s_{\mu}-s_{\nu}}{s_{\mu}-1}
\end{equation}
provided that $s_{\mu}>1$. Note that since we are in the regime $1 \leq p \leq 2$, in order that the interval of convergence of the series be not empty, we need $\frac{d+s_{\mu}-s_{\nu}}{s_{\mu}-1} < 2$, which is equivalent to $s_{\mu}+s_{\nu}>d+2$.

\vskip.125in 

Now suppose that $2 \leq p<\infty$. We now use the local smoothing estimates from \cite{MSS93}. The estimate (\ref{estI}) now takes the form 
$$ I \leq C 2^{-j \frac{d-1}{2}} \cdot 2^{j (d-1)\left| \frac{1}{2}-\frac{1}{p}\right|} 2^{-j \epsilon(p)} 2^{\frac{j}{p}} 2^{j \delta}
{||P_j(f\mu)||}_{L^p({\Bbb R}^d)}$$
$$ \leq C_\eta 2^{-j \frac{d-1}{p}} 2^{-j \epsilon(p)} 2^{\frac{j}{p}} 2^{j \delta} 2^{j \frac{d-s_{\mu} + \eta}{p'}}{||f||}_{L^p(\mu)}=C_\eta 2^{-j \frac{d-2}{p}} 2^{-j \epsilon(p)} 2^{j \delta} 2^{j \frac{d-s_{\mu} + \eta}{p'}}{||f||}_{L^p(\mu)}. $$

Plugging in (\ref{Lpprime}) we obtain 
$$ I \cdot II \leq C_\eta2^{-j \frac{d-2}{p}} 2^{-j \epsilon(p)} 2^{j \frac{d-s_{\mu}}{p'}} 2^{j (\delta + \eta)} 2^{j \frac{d-s_{\nu}}{p}}{||f||}_{L^p(\mu)} 
{||g||}_{L^{p'}(\nu \times \gamma)}.$$

Recalling that $\delta$ and $\eta$ are arbitrarily small once again, we see that the geometric series converges if 
$$ d-\frac{s_{\mu}}{p'}-\frac{s_{\nu}}{p}<\frac{d-2}{p}+\epsilon(p),$$ which reduces to 
\begin{equation}\label{ConditionForConvergenceIfPBiggerThan2}
 d-s_{\mu}+\frac{s_{\mu}-s_{\nu}}{p}<\frac{d-2}{p}+\epsilon(p).
\end{equation}

\vskip.125in 

\subsubsection{Proof of parts ii) and iii)} 

It is just a matter of checking that the convergence conditions \eqref{ConditionForConvergenceIfPLessThan2} when $1 \leq p \leq 2$, and \eqref{ConditionForConvergenceIfPBiggerThan2} when $2 \leq p < \infty$ hold.

\qed

\vskip.25in 

\subsection{Proof of Theorem \ref{steinexample} and Theorem \ref{mattilaexample} } 

\subsubsection{Proof of Theorem \ref{steinexample}} 
First assume $s>0$. Let $s_{\mu} \equiv s$, $B = B(0,1)$ be the unit ball in $\mathbb{R}^d$, and define $d\mu(x)={|x|}^{-d+s}dx$, $t=|x|$ and 
$$f(x)={|x|}^{1-s} \log^{-1} \left( \frac{1}{|x|} \right) \chi_{\frac{1}{2}B}(x).$$ 

Then 

$$ \int {|x|}^{p(1-s)}  \log^{-p} \left( \frac{1}{|x|} \right) \chi_{\frac{1}{2}B}(x) d\mu(x)=\int |x|^{p(1-s)-d+s}  \log^{-p} \left( \frac{1}{|x|} \right) \chi_{\frac{1}{2}B}(x)dx<\infty$$ if $1 < p \leq \frac{s}{s-1}$ in case $s>1$, and for all $1 \leq p < \infty$ if $0 < s \leq 1$. 

\vskip.125in

On the other hand, for $x \in \frac{1}{4}B$,

$$ A_{|x|}(f\mu)(x)=\int { \left( |x-|x|y| \right) }^{1-s} {(|x-|x|y|)}^{-d+s} \log^{-1} \left( \frac{1}{|x-|x|y|} \right) d\sigma(y)$$
$$=\int { \left( |x-|x|y|  \right) }^{-d+1} \log^{-1} \left( \frac{1}{|x-|x|y|} \right) d\sigma(y) \equiv \infty$$

since the sphere is a $d-1$-dimensional surface. After a simple renormalization (or redefining the supremum in the operator  ${\mathcal A}(f\mu)$ in Theorem \ref{maximalmain} to be over $t \in [\frac{1}{8}, \frac{1}{4}]$, say), we conclude that ${\mathcal A}(f\mu)$ is, in general, infinite, if $f \in L^p(\mu)$ for $1 < p \leq \frac{s}{s-1}$ in case $s > 1$, and for $1 \leq p < \infty$ in case $0 < s \leq 1$. 

Now assume $s=0$.  Let $d\mu(x)={|x|}^{-d}\log^{-u}\left(\frac{1}{|x|}\right)dx$, $t=|x|$ and 
$$f(x)={|x|} \log^{\beta} \left( \frac{1}{|x|} \right) \chi_{\frac{1}{2}B}(x),$$ 
for some $u = \beta >1$. Then the measure $\mu$ satisfies \eqref{polygrowth} with $s_{\mu} = 0$ and the same argument as the case $s>0$ gives that $f \in L^p (\mu)$ for any $1 < p < \infty$, yet ${\mathcal A}(f\mu)(x)= \infty$.

\vskip.125in 

\subsubsection{Proof of Theorem \ref{mattilaexample}} 

\vskip.125in 

Consider $E=C_{\alpha}^{d-1} \times C_{\beta}$, where $C_{\alpha} \subset [0,1]$ is a Cantor-type set of dimension $\alpha$ with gauge function $h(t) = t^\alpha$, if $0 < \alpha <1$. If $\alpha = 1$ we take $C_{\alpha} = [0,1]$, and if $\alpha = 0$ we either take $C_{\alpha}$ with $\alpha \approx 0$ and let $\alpha \to 0^{+}$ in the arguments below, or else take a Cantor set of dimension $0$ with gauge function $h(t)$ and adapt the arguments below in the obvious way to that gauge function (see e.g. \cite{Mat95}, sects. 4.9 and 4.11). Let $\mu = \mu_\alpha \times \dots \times \mu_\alpha \times \mu_\beta$ be a Frostman measure for $E$, where $\mu_\alpha$ and $\mu_\beta$ are Frostman measures for  $C_{\alpha}$ and $ C_{\beta}$ respectively.  Let 
$f(y)={|y_d|}^{-\frac{\beta}{p}}$. Recall that a $\sqrt{\epsilon} \times \dots \times \sqrt{\epsilon} \times \epsilon$ rectangle fits inside an annulus of radius $1$ and width $\approx \epsilon$. After setting $t=x_d$, we see that for $x \in E$,

$$ {\mathcal A}f(x) \gtrsim \liminf_{\epsilon \to 0^{+}}  \epsilon^{-1} \int_{x_d - \epsilon \leq |x-y| \leq x_d+\epsilon} f(y) d\mu(y) \gtrsim 
\liminf_{\epsilon \to 0^{+}} \epsilon^{-1} \epsilon^{-\frac{\beta}{p}} \epsilon^{\frac{\alpha(d-1)}{2}+\beta}$$
\begin{equation} \label{cantorstream} =\liminf_{\epsilon \to 0^{+}} \epsilon^{-1+\frac{\alpha(d-1)}{2}+\frac{\beta}{2}+\frac{\beta}{2}-\frac{\beta}{p}}= 
\liminf_{\epsilon \to 0^{+}} \epsilon^{-1+\frac{s}{2}+\beta \left( \frac{1}{2}-\frac{1}{p} \right)}. \end{equation} 

If $1 \leq p \leq 2$, $\frac{1}{2}-\frac{1}{p}\leq 0$, so we set $\beta=1$ and obtain from (\ref{cantorstream}): 
$$ \epsilon^{-1+\frac{s}{2}+\frac{1}{2}-\frac{1}{p}}=\epsilon^{\frac{s}{2}-\frac{1}{2}-\frac{1}{p}}$$ and part a) follows by taking $\mu = \nu$. 

To prove part b) note that we are in the range $p \ge 2$, so $\frac{1}{2}-\frac{1}{p} \ge 0$. Let $\alpha=1$ and $\beta=s-1$. Plugging this into (\ref{cantorstream}) yields the conclusion again by taking $\mu = \nu$. 

To prove part c) just take $\beta=0$ in (\ref{cantorstream}) and again take $\mu = \nu$.


\vskip.25in 

\section{Applications to the wave equation} 

\vskip.125in 

In this section we work out some applications of the results from the previous section to the wave equation in three dimensions. Our results can be easily extended to other dimensions as well, but we mostly stick to the three dimensional setup for the sake of ease of presentation. We consider the initial value problem (\ref{wave}): 
$$ \Delta u=\frac{\partial^2 u}{\partial t^2}; \ u(x,0)=0; \ \frac{\partial u}{\partial t}(x,0)=f(x).$$

As we note in the introduction, 
$$ u(x,t)=ct \int_{S^2} f(x-ty) d\sigma(y),$$ where $\sigma$ is the Lebesgue measure on the sphere, as before. 

Suppose that we consider a slightly modified version of the same initial value problem 
\begin{equation} \label{waveweirdproblem} \Delta u=\frac{\partial^2 u}{\partial t^2}; \ u(x,0)=0; \ \frac{\partial u}{\partial t}(x,0)=f\mu(x), \end{equation} where $\mu$ is a compactly supported Borel measure. With a similar proof as in the previous section (but since we consider a fixed $t$ there is no Sobolev embedding or local smoothing, which makes the calculations simpler), we get the following Theorem.

\begin{theorem} \label{waveweird} Consider the initial value problem (\ref{waveweirdproblem}) for $x \in \mathbb{R}^3$ where $\mu$ is a compactly supported Borel measure satisfying $\mu(B(x,r)) \leq Cr^{s_{\mu}}$ for some $s_{\mu}>0$. Suppose that $\nu$  is a compactly supported Borel measure satisfying $\nu(B(x,r)) \leq Cr^{s_{\nu}}$ for some $s_{\nu}>0$. Then for every $t>0$, there exists $C=C(t)>0$ such that 
\begin{equation} \label{L2wave} {||u(\cdot, t)||}_{L^2(\nu)} \leq C{||f||}_{L^2(\mu)} \ \text{if} \ s_{\mu}+s_{\nu}>4. \end{equation} 

When $p \ge 2$, 
\begin{equation} \label{Lp+wave} {||u(\cdot, t)||}_{L^p(\nu)} \leq C{||f||}_{L^p(\mu)} \ \text{if} \ \frac{s_{\mu}}{p'}+\frac{s_{\nu}}{p}>3-\frac{2}{p}. \end{equation} 

When $p \leq 2$, 
\begin{equation} \label{Lp-wave} {||u(\cdot, t)||}_{L^p(\nu)} \leq C{||f||}_{L^p(\mu)} \ \text{if} \ \frac{s_{\mu}}{p'}+\frac{s_{\nu}}{p}>1+\frac{2}{p}. \end{equation} 

\end{theorem} 

The $L^p$ estimates in (\ref{Lp+wave}) and (\ref{Lp-wave}) are interesting in their own right, so we state a higher dimensional analog. 

\begin{theorem} \label{avoplp} Let $A_tf(x)=\sigma_t*(f\mu)(x)$, where, as before, $\sigma_t$ is the surface measure on the sphere of radius $t$ in ${\Bbb R}^d$, $d \ge 2$, $\mu$ is a compactly supported Borel measure such that $\mu(B(x,r)) \leq Cr^{s_{\mu}}$ and $\nu$ is a compactly supported Borel measure such that $\nu(B(x,r)) \leq Cr^{s_{\nu}}$. 

\vskip.125in 

When $p \ge 2$, then for every fixed $t>0$, 
\begin{equation} \label{Lp+av} {||A_tf||}_{L^p(\nu)} \leq C(t){||f||}_{L^p(\mu)} \ \text{if} \ \frac{s_{\mu}}{p'}+\frac{s_{\nu}}{p}>d-\frac{d-1}{p}. \end{equation} 

When $p \leq 2$, then for every fixed $t>0$, 
\begin{equation} \label{Lp-av} {||A_tf||}_{L^p(\nu)} \leq C(t){||f||}_{L^p(\mu)} \ \text{if} \ \frac{s_{\mu}}{p'}+\frac{s_{\nu}}{p}>1+\frac{d-1}{p}. \end{equation}  \end{theorem} 

\vskip.125in 

\begin{remark} In particular, in the same way as with Corollary \ref{BlowUpSetforThmMain}, by taking $\mu$ to be Lebesgue measure in $\mathbb{R}^d$ restricted to an appropriate large compact set (say a ball of radius $> 4t$), and taking $\nu$ to be first a restriction of Lebesgue measure to a compact set, and then a Frostman measure on an arbitrary set of Hausdorff dimension $>1$, we see that Theorem \ref{waveweird} implies that  $u(x,t)$ is an $L^2({\Bbb R}^d)$ function in the $x$ variable that is well-defined up to a set of Hausdorff dimension $=1$, for every fixed $t>0$. \end{remark} 

\vskip.125in 

\begin{corollary} \label{avoplpcor} [Blow-up sets for the spherical averaging operator] Let $A_tf(x)$ be as above and define 
$$ E_f=\{x \in {\Bbb R}^d: A_tf(x)=\infty \}.$$ 

If $f \in L^p({\Bbb R}^d)$, $1 \leq p \leq 2$, then 
\begin{equation} \label{<2cor} dim_{{\mathcal H}}(E_f) \leq d-(p-1)(d-1). \end{equation} 

If $f \in L^p({\Bbb R}^d)$, $p \ge 2$, then 
$$ dim_{{\mathcal H}}(E_f) \leq 1.$$  

\end{corollary} 

It is not difficult to see that Corollary \ref{avoplpcor} follows from Theorem \ref{avoplp} in the same way as Corollary \ref{BlowUpSetforThmMain} follows from Theorem \ref{main}. Let us now consider the extent to which Corollary \ref{avoplpcor} is sharp. Let 
$$ f(x)={|x|}^{-(d-1)} \log^{-1} \left( \frac{1}{|x|} \right) \chi_{\frac{1}{2}B}(x),$$ where $B$ is the unit ball. Then $f \in L^p(B)$ for $p \leq \frac{d}{d-1}$. On the other hand, 
$$ A_1f(x) \approx \int_{S^{d-1}} {|x-y|}^{-(d-1)} \log^{-1} \left( \frac{1}{|x-y|} \right) d\sigma(y) \equiv \infty \ \text{for} \ x \in S^{d-1}.$$ 

It follows that $dim_{{\mathcal H}}(E_f)=d-1$, which matches Corollary \ref{avoplpcor} since plugging $p=\frac{d}{d-1}$ into $(\ref{<2cor})$ yields $dim_{{\mathcal H}}(E_f) \leq d-1$.

\vskip.125in

\begin{remark}
We could also consider $ {\mathcal A}f(x)$ and corresponding estimates for the blowup set in $(x,t) \in \mathbb{R}^d \times \mathbb{R}^{+}$. Then, if one assumes the local smoothing conjecture, one gets better estimates than what one gets with the known local smoothing estimates used above. Conversely, these estimates yield a possible strategy for disproving the local smoothing conjecture, by finding appropriately large blowup sets. We shall explore this issue in more detail in the sequel. \end{remark}

\vskip.125in 


\end{document}